\documentclass[preprint,12pt]{elsarticle}

\usepackage{amsmath,amssymb,amscd}
\usepackage{amsthm}
\usepackage{subfigure}
\usepackage[nodots]{numcompress}
\usepackage{verbatim}
\usepackage{fancybox}
\usepackage{subfigure}
\usepackage[usenames,dvipsnames]{color}
\usepackage{sidecap}
\usepackage{caption}
\usepackage{ifthen}
\usepackage{textcomp}
\usepackage[all]{xy}

\newcommand{\R}{{\mathbb R}}

\newcommand{\Q}{{\mathbb Q}}

\newcommand{\dem}{{\em Proof: \;}}
\newcommand{\fdem}{\hfill $\square$}

\theoremstyle{plain}
\newtheorem{teo}{Theorem}[section]
\newtheorem{lema}[teo]{Lemma}
\newtheorem{cor}[teo]{Corollary}
\newtheorem{prop}[teo]{Proposition}

\theoremstyle{definition}
\newtheorem{defi}{Definition}[section]

\theoremstyle{remark}

\begin{document}

\begin{frontmatter}

\title{Factorization of network reliability with perfect nodes I: Introduction and Statements}

\author[IMATE]{Juan Manuel Burgos}
\address[IMATE]{Instituto de Matem\'aticas, Universidad Nacional Aut\'onoma de M\'exico, Unidad Cuernavaca. Av. Universidad s/n, Col. Lomas de
Chamilpa. Cuernavaca, Morelos M\'exico, 62209.\\ \texttt{\scriptsize{Email: burgos@matcuer.unam.mx}}}

\author[lpe]{Franco Robledo Amoza}
\address[lpe]{
     Laboratorio de Probabilidad y Estad\'{\i}stica\\ Facultad de Ingenier\'ia, Universidad de la Rep\'ublica\\
     Julio Herrera y Reissig 565\\ 11300, Montevideo, Uruguay.\\ \texttt{\scriptsize{Email: frobledo@fing.edu.uy}}}

\begin{abstract}
A new general all terminal network reliability factorization theorem is stated. We relegate the proof to a forthcoming second part paper.
\end{abstract}


\end{frontmatter}

%-------------------------------------------------------------------------

\section{Introduction}

The network reliability factorization theorem [Mo] and the reduction transformations (series-parallel, polygon-to-chain [Wo] and delta-star [Ga]) are the key stone of the known factoring algorithms [SC] for the exact calculation of network reliability. Beside these results and the well known factorization through an articulation point, no other general factorization theorem is known in exact network reliability calculation. This paper gives a new general all terminal network reliability factorization theorem solving the following problem:

\textbf{Problem:} \emph{Given a decomposition of a stochastic graph $G$ by subgraphs $G_{1}$ and $G_{2}$ only sharing nodes, express the reliability of $G$ in terms of the reliabilities of the graphs resulting from $G_{1}$ and $G_{2}$ identifying the common nodes shared by them in all possible ways.}

%Exact calculation of $K$-network reliability is a NP-hard problem [Co] so it is necessary to know how to do parallel processing in order to reduce the processing time. This theorem has the following immediate parallel processing application: Suppose we are willing to calculate the exact network reliability of a stochastic graph $G=(V,E,(p_{e})_{e\in E})$ with $m$ edges and perfect nodes. A priori the processing time involve $2^{m}$ steps. Cutting the graph $G$ through a set of nodes obtaining $G_{1}$ and $G_{2}$ with approximately $m/2$ edges each, the solution of the above problem allow us to calculate the reliability of $G$ with $2^{m/2}$ steps in each processing line.

\begin{figure}
\begin{center}
  \includegraphics[width=1.0\textwidth]{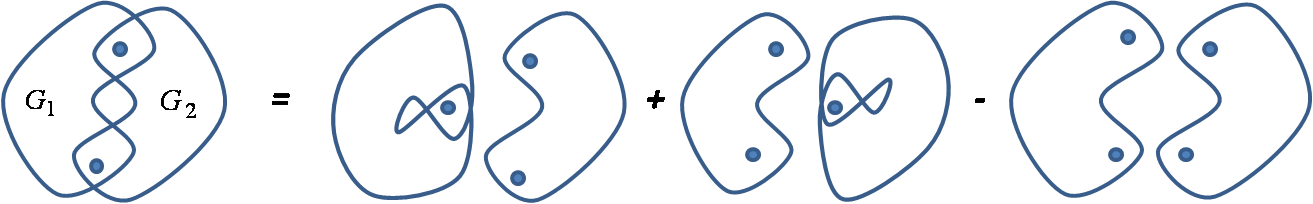}\\
  \end{center}
  \caption{}\label{nDos}
\end{figure}

As an example, consider a graph $G$ with the property it that can be decomposed by subgraphs $G_{1}$ and $G_{2}$ only sharing a pair nodes $a$ and $b$. Denoting $\hat{G}_{i}$ the graph resulting from the identification of the nodes $a$ and $b$ in the graph $G_{i}$, we have the following factorization formula for the all terminal reliability: $$R(G)=R(G_{1})R(\hat{G}_{2}) + R(\hat{G}_{1})R(G_{2}) - R(G_{1})R(G_{2})$$
This formula is illustrated in Figure \ref{nDos} (see Figure \ref{nTres} for the case of three sharing nodes). A simple proof of this particular case is given in appendix A. After the preliminary section (second section), the rest of the paper is devoted to the formulation of the general case (third and fourth sections). We relegate the proof of the theorem to the second part of this paper [Bu]: Factorization of network reliability with perfect nodes II.

We remark that the paper is about a new purely mathematical theorem and it has nothing to do with computational complexity nor algorithms. Mathematically, the proof introduces novel techniques based on an algebraic approach and connectivity properties.

\section{Preliminaries}

%Def de grafo estoc'astico
%Def de estados, Pathset, Cutset.
%Def de Confiabilidad
%Def. de contracci'on y eliminaci'on de una arista
%F'ormula de fact.simple
%Ejemplo

\subsection{Reliability}

The mathematical model of a network whose nodes are perfect and its edges can fail is a stochastic graph [Co], [Sh]; i.e. an undirected graph with associated Bernoulli variables to its edges. I assume throughout the paper that the forward slash means ``such that"`.

\begin{defi}
An undirected graph $G$ is $(V,E,\phi)$ such that $V$ and $E$ are are finite sets whose elements will be called nodes and edges respectively $E$ and $\phi$ is a function from $E$ to $\{\{a,b\}\ /\ a,b\in V\}$ specifying the nodes attached to each edge. We will make an abuse of notation and denote the graph $G$ simply as $G=(V,E)$. Nodes and edges of $G$ will be denoted by $V(G)$ and $E(G)$ respectively.
\end{defi}

%\begin{defi}
%The graph $G'=(V',E',\phi')$ is a subgraph of $G=(V,E,\phi)$ if $v'\subset V$, $E'\subset E$ and $\phi'= \phi|_{E'}$.
%\end{defi}

\begin{defi}
A stochastic graph $G$ is $(V,E,\phi, \Phi)$ such that $(V, E,\phi)$ is a graph and $\Phi:E\rightarrow Ber$ is a function which associates a Bernoulli variable to each edge in such a way that these variables are independent.
\end{defi}

\begin{defi}
The graph $G'=(V',E',\phi', \Phi')$ is a subgraph of $G=(V,E,\phi, \Phi)$ if $V'\subset V$, $E'\subset E$, $\phi'= \phi|_{E'}$ and $\Phi'= \Phi|_{E'}$.
\end{defi}

Each Bernoulli variable is characterized by a parameter $p$ in the $[0,1]$ closed interval and we can write a stochastic graph as $(G, \{p_{e}\}_{e\in E})$ where $G$ is an undirected graph and $p_{e}$ is the parameter of the variable $\Phi(e)$.

\begin{defi}
A state $\mathcal{E}$ of the graph $G=(V,E)$ is a function $\mathcal{E}: E\rightarrow \{0, 1\}$. An edge $e$ will be called operative if $\mathcal{E}(e)=1$ and will be called non-operative otherwise.
\end{defi}

A state $\mathcal{E}$ of the graph $G$ will be called a PathSet (or operative) if the subgraph $G_{\mathcal{E}}$ resulting from the removal of the non-operative edges of $G$ is spanning and edge connected. Otherwise the state will be called a CutSet (or non operative).

\begin{defi}\label{DefConfiabilidad}
The Reliability of a stochastic graph $G$ is $$R(G)=P(\mathcal{E}\ is\ a\ PathSet)$$
\end{defi}

Because of the independence of the Bernoulli variables associated to the edges, we can calculate the Reliability in the following way:

\begin{equation}\label{ProbabilidadPathset}
    P(\mathcal{E})=\prod_{e\in E(G)}p_{e}^{\mathcal{E}(e)}(1-p_{e})^{1-\mathcal{E}(e)}
\end{equation}

$$R(G)=\sum_{\mathcal{E}\ is\ a\ PathSet} P(\mathcal{E})$$

\subsection{Simple Factorization}

\begin{figure}
\begin{center}
  \includegraphics[width=0.4\textwidth]{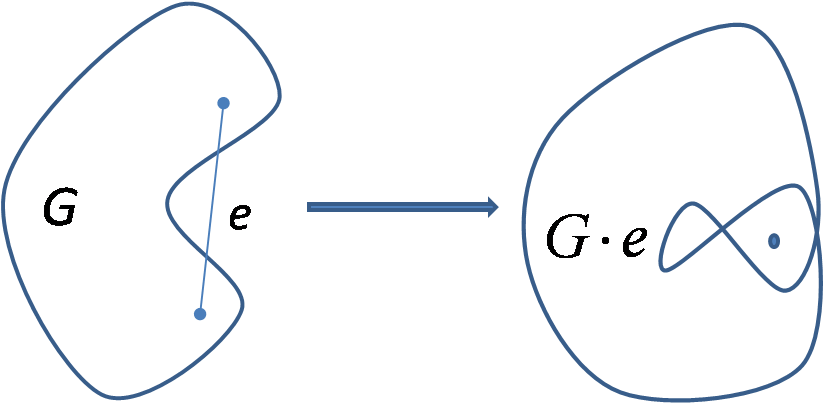}\\
  \end{center}
  \caption{Contraction of the edge \textit{e}}\label{contrac}
\end{figure}

\begin{defi}
Consider an edge $e\in E$ of a graph $G=(V,E, \phi)$ and define the following equivalence relation in $V$: $a\sim b$ if $a=b$ or $\{a,b\}=\phi(e)$. Consider the suryective canonical function $\pi:V\rightarrow V/\sim$ such that $\pi(a)=[a]_{\sim}$. We define the contraction of an edge $e$ in $G$ as the graph $G\cdot e$ such that $G\cdot e = (V/\sim, E-\{e\}, \bar{\phi})$ where $\bar{\phi}(e)= \{\pi(a),\pi(b)\}$ if $\phi(e)= \{a,b\}$ (see Figure \ref{contrac}).
\end{defi}

\begin{figure}
\begin{center}
  \includegraphics[width=0.4\textwidth]{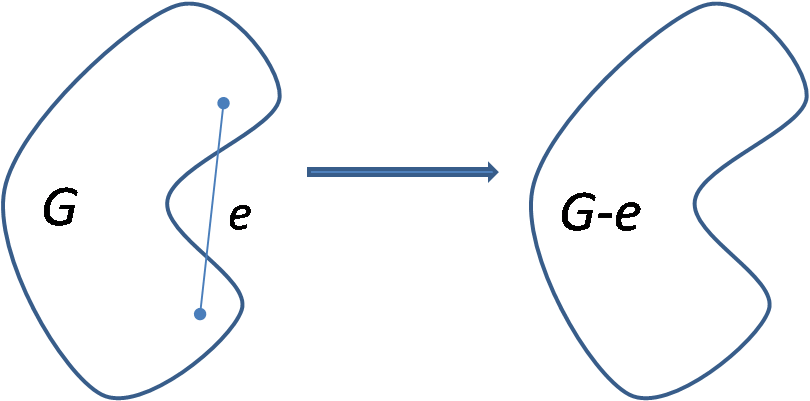}\\
  \end{center}
  \caption{Deletion of the edge \textit{e}}\label{extrac}
\end{figure}

\begin{defi}
Consider an edge $e\in E$ of the graph $G=(V,E)$. We define the deletion of the edge $e$ of $G$ as the graph (see Figure \ref{extrac}) $G-e=(V,E-\{e\}, \phi )$
\end{defi}

The following is the simple factorization proposed for the first time by Moskovitz [Mo]:

\begin{prop}\label{FactSimple}
Consider a stochastic graph $G$. For any edge $e$ of $G$ we have that $$R(G)=p_{e}R(G\cdot e) + (1-p_{e})R(G-e)$$
\end{prop}

\begin{defi}
Define the following partial order in the set of states of $G$: $\mathcal{E} \leq \mathcal{F}$ if $\mathcal{E}^{-1}(1) \subset \mathcal{F}^{-1}(1)$. The state $\mathcal{E}$ is a minpath if it is minimal in the set of PathSets.
\end{defi}

The PathSets have the following coherence property: If $\mathcal{E} \leq \mathcal{F}$ and $\mathcal{E}$ is a PathSet, then $\mathcal{F}$ is a PathSet. For every PathSet $\mathcal{F}$ there is a minpath $\mathcal{E}$ such that $\mathcal{E} \leq \mathcal{F}$. We conclude that the set of PathSets equals the set of states $\mathcal{F}$ greater or equal to some minpath. This motivates the following definition [Co]:

\begin{defi}
An edge $e$ of a graph $G$ will be called irrelevant if $\mathcal{E}(e)=0$ for each minpath $\mathcal{E}$.
\end{defi}

\begin{prop}
If $e$ is an irrelevant edge of a stochastic graph $G$, then $$R(G)=R(G-e)$$
\end{prop}

%\subsection{Reliability Polynomial}

Consider a graph $G$ and the stochastic graph $G_{p}$ whose underlying graph is $G$ and its Bernoulli variables are identical and independent with parameter $p$. By proposition \ref{FactSimple} and the finiteness of the graph, we have that the reliability of $G$ is a polynomial in the parameter $p$:

\begin{defi}
$R(p)= R(G_{p})$ is the reliability polynomial of the graph $G$; i.e. $$R(p)=\sum_{i=0}^{m}C_{i}p^{i}(1-p)^{m-i}$$ where $C_{i}$ is the number of PathSets with exactly $i$ operative edges where $m$ is the number of edges in $G$.
\end{defi}

The following complexity theorem is due to Provan and Ball [PB]:
%In particular, the calculation of the reliability polynomial is equivalent to the calculation of the numbers $C_{i}$. Each minpath is a Steiner spanning tree of $G$ and the calculation of the reliability polynomial solves the following $NP$-complete problem: Given a number $b$, Is there a Steiner spanning tree with edge number less than or equal to $b$? To answer this question we just have to see weather or not $m\leq b$ where $C_{m}$ is the first non zero coefficient of the above expression. Then, the calculation of the reliability polynomial is at least $NP$-complete. In fact, we have the following [Co]:

\begin{teo}
The calculation of the reliability polynomial is a $NP$-hard problem.
\end{teo}

\section{Combinatorics of the Problem}

Consider a stochastic graph $G= (V,E,\phi, \Phi)$ such that there exist a pair of subgraphs $G_{1}$ and $G_{2}$ of $G$ such that $E=E_{1}\cup E_{2}$, $V=V_{1}\cup V_{2}$, $E_{1}\cap E_{2}= \emptyset$ and $V_{1}\cap V_{2}= \{v_{1}, v_{2},\ldots v_{n}\}$. Because $E_{1}$ and $E_{2}$ is a disjoint union of edges, we can decompose the set of states of $G$ in terms of the states of $G_{1}$ and $G_{2}$ as follows: $$\{0, 1\}^{E}= \{0, 1\}^{E_{1}\cup E_{2}}= \{0, 1\}^{E_{1}}\times \{0, 1\}^{E_{2}}$$ and we can write a state $\mathcal{E}$ as a pair $\mathcal{E}= (\mathcal{E}_{1}, \mathcal{E}_{2})$ where $\mathcal{E}_{i}\in \{0, 1\}^{E_{i}}$, $i=1,2$. According to this decomposition we have $$P(\mathcal{E})=P(\mathcal{E}_{1})\ P(\mathcal{E}_{2})$$

For every state of $G$ we associate a partition of $\{v_{1}, v_{2},\ldots v_{n}\}$ as follows: Consider the equivalence relation: $v_{i}\sim v_{j}$ if there is a path in  $G_{\mathcal{E}}$ joining $v_{i}$ with $v_{j}$. We define the partition $$\left[\mathcal{E}\right]=\{v_{1},v_{2},\ldots v_{n}\}/\sim$$ and we call it the partition of $\mathcal{E}$ with respect to $G$. Analogously, we define the partition of $\mathcal{E}_{i}$ with respect to the subgraph $G_{i}$ and denote it by $\left[\mathcal{E}_{i} \right] $, $i=1,2$. We will call the partition $\{\{v_{1},v_{2},\ldots v_{n}\}\}$ the trivial partition.

\begin{figure}
\begin{center}
  \includegraphics[width=0.6\textwidth]{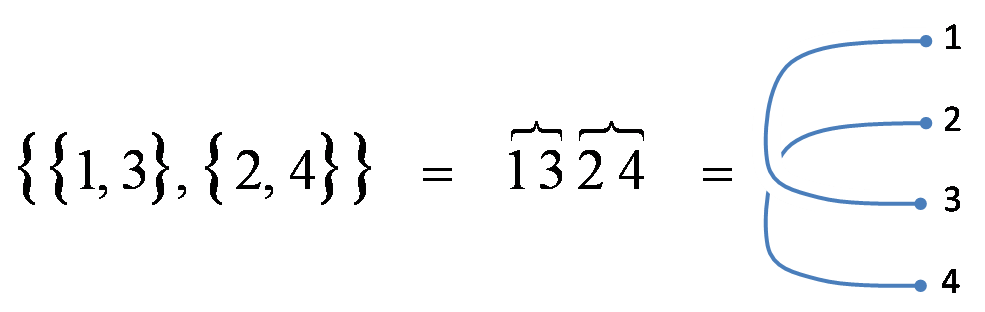}\\
  \end{center}
  \caption{Connectivity State}\label{EjEstadoConect}
\end{figure}

Denote by $Part_{n}$ the set of partitions of $\{v_{1},v_{2},\ldots v_{n}\}$. Figure \ref{EjEstadoConect} shows some useful notational and diagrammatical  ways to represent a partition. The set $Part_{n}$ has a monoid structure with unit $\{ \{v_{1} \} , \{ v_{2} \} , \ldots \{ v_{n} \} \}$ under the following product: given the partitions $\mathcal{A}$ and $\mathcal{B}$, the product $\mathcal{A}\cdot \mathcal{B}$ is the finer partition that is coarser to $\mathcal{A}$ and coarser to $\mathcal{B}$. Observe that the product of any partition with the trivial one is trivial.

Is clear that if the partition of $\mathcal{E}$ with respect to $G$ is not the trivial one then $P(\mathcal{E})=0$ and the state $\mathcal{E}$ doesn't contribute in the reliability of $G$. This motivates the following definition:

\begin{defi}
The connectivity state of a state $\mathcal{E}$ of $G$ is the pair of partitions $\left(\left[\mathcal{E}_{1} \right],\ \left[\mathcal{E}_{2} \right] \right)$ and we say it is connected if the partition $\left[\mathcal{E}\right]$ is the trivial one.
\end{defi}

\begin{figure}
\begin{center}
  \includegraphics[width=0.4\textwidth]{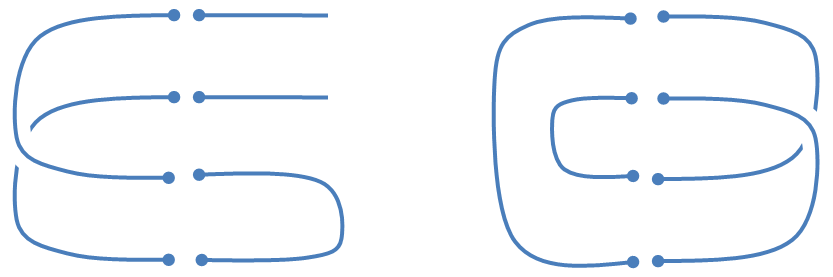}\\
  \end{center}
  \caption{Connected connectivity states of $G$}\label{EjEstadoConexo}
\end{figure}

\begin{figure}
\begin{center}
  \includegraphics[width=0.4\textwidth]{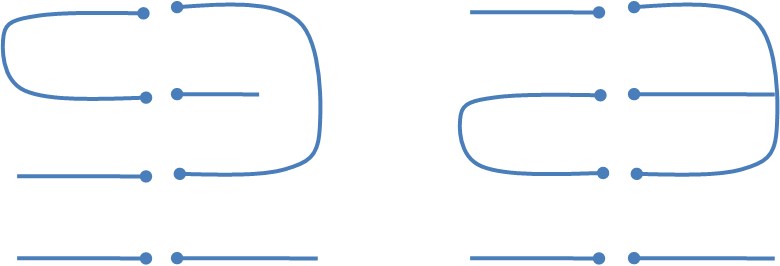}\\
  \end{center}
  \caption{Non connected connectivity states of $G$}\label{EjEstadoNoConexo}
\end{figure}

Figures \ref{EjEstadoConexo} and \ref{EjEstadoNoConexo} show examples of connectivity states of $G$. This way we have the following formula for the reliability of $G$: $$R(G)=\sum_{(\left[\mathcal{E}_{1}\right], \left[\mathcal{E}_{2}\right])\ connected} P(\mathcal{E})= \sum_{(\left[\mathcal{E}_{1}\right], \left[\mathcal{E}_{2}\right])\ connected} P(\mathcal{E}_{1})\ P(\mathcal{E}_{2})$$

We can rewrite the above identity only in terms of the partitions as follows: Given a partition $\mathcal{A}$ of $\{v_{1},v_{2},\ldots v_{n}\}$, define $$P_{i}(\mathcal{A})= \sum_{\left[\mathcal{E}_{i}\right]= \mathcal{A}} P_{i}(\mathcal{E}_{i})$$ This way $$R(G)=\sum_{(\mathcal{A}, \mathcal{B})\ connected} P_{1}(\mathcal{A})P_{2}(\mathcal{B})$$ where $\mathcal{A}$ and $\mathcal{B}$ are partitions of $\{v_{1},v_{2},\ldots v_{n}\}$ and we say the pair $(\mathcal{A}, \mathcal{B})$ is connected if $\mathcal{A}\cdot \mathcal{B}$ is the trivial partition; i.e. the resulting diagram is connected (see Figures \ref{EjEstadoConexo} and \ref{EjEstadoNoConexo}). We have written the reliability of $G$ in terms of the partitions only and not in terms of particular states. 

The above formula was also given in [Ro] (under a different notational scheme) where an algorithm for the reliability exact calculation is developed based on it. The remarkable property of the above formula is that now we can express the reliability of $G$ in terms of reliabilities of particular states of the subgraphs, something impossible without taking off the non connected connectivity states in the sum. However, the above formula is not a factorization theorem in the sense that it doesn't express the reliability of $G$ in terms of the reliabilities of the subgraphs, our original problem. In order to do so, we need some machinery first. The rest of the section is devoted to the translation of our original probabilistic problem given in the introduction to a purely combinatorial one.

\begin{defi}
For each partition $\mathcal{A}$ of $\{v_{1}, v_{2},\ldots v_{n}\}$ denote by $G^{\mathcal{A}}_{i}$ the graph resulting from the identification of the nodes in $\{v_{1},v_{2},\ldots v_{n}\}$ of $G_{i}$ by the partition $\mathcal{A}$.
\end{defi}

This formula suggests the following algebraic construction: Consider the $\Q$-vector space $V$ generated by $Part_{n}$ and the linear functional $P_{i}:V\rightarrow \R$ which is the linear extension of the probability $P_{i}$ (we are abusing of the notation and denoting the functional by the same name). This way we get the functional $$P_{1}\otimes P_{2}:\ V\otimes V \rightarrow \R$$ and the following expression for the reliability: $$R(G)=P_{1}\otimes P_{2}\left(\sum_{(\mathcal{A}, \mathcal{B})\ connected} \mathcal{A} \otimes \mathcal{B}\right)$$

\begin{defi}
Denote by $R_{\mathcal{A}}(G_{i})$ the reliability of $G^{\mathcal{A}}_{i}$
\end{defi}

\begin{lema}
$R_{\mathcal{A}}(G_{i})=\sum_{\mathcal{B}\ /\ (\mathcal{A}, \mathcal{B})\ connected} P_{i}(\mathcal{B})$
\end{lema}
\dem
The set of states and Bernoulli variables of $G_{i}$ and $G^{\mathcal{A}}_{i}$ are equal. Because of the definition of reliability, the result follows from the following observation: A state $\mathcal{E}_{i}$ of $G_{i}$ is a PathSet of $G^{\mathcal{A}}_{i}$ if and only if $(\mathcal{A}, \left[\mathcal{E}_{i}\right])$ is connected.
\fdem

\begin{figure}
\begin{center}
  \includegraphics[width=0.8\textwidth]{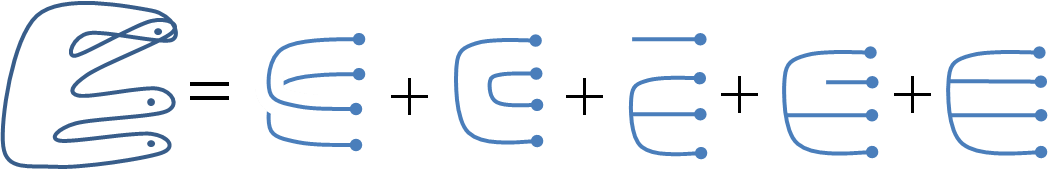}\\
  \end{center}
  \caption{}\label{EjRunodos}
\end{figure}

We can rewrite the above formula with the functional $P_{i}$: $$R_{\mathcal{A}}(G_{i})=P_{i}\left(\sum_{\mathcal{B}\ /\ (\mathcal{A}, \mathcal{B})\ connected} \mathcal{B}\right)$$ Then we have finally translated the original probabilistic problem in the following purely combinatorial one:

\textbf{Problem (Combinatorial Version):} \emph{Express the vector $R\in V\otimes V$ given by $$R= \sum_{(\mathcal{A}, \mathcal{B})\ connected} \mathcal{A} \otimes \mathcal{B}$$ in terms of the vectors $R_{\mathcal{A}}\in V$ given by (see Figure \ref{EjRunodos}) $$R_{\mathcal{A}}=\sum_{\mathcal{B}\ /\ (\mathcal{A}, \mathcal{B})\ connected} \mathcal{B}$$
}

%Hacer otra figura con la identificaci'on de v'ertices

The relation between the combinatorial and the probabilistic problem is given by the functional $P_{1}\otimes P_{2}:\ V\otimes V \rightarrow \R$; i.e. Acting on the combinatorial expression by $P_{1}\otimes P_{2}$ gives the probabilistic one.

\section{The Factorization Theorem}

Having translated our problem into a combinatorial one, we will develop its solution in this combinatorial framework. Considering an ordering in the basis $Part_{n}$ we can write $$R= \sum_{(\mathcal{A}, \mathcal{B})\ connected} \mathcal{A} \otimes \mathcal{B} =\sum_{i,j=1}^{m} a_{ij}\mathcal{A}_{i}\otimes \mathcal{A}_{j}$$ where $m$ is the cardinality of $Part_{n}$ and $A=(a_{ij})$ is the matrix given by $a_{ij}=1$ if $(\mathcal{A}_{i}, \mathcal{A}_{j})$ is connected and $a_{ij}=0$ if it is not. The matrix $A$ will be called the \emph{connectivity matrix}. Consider the linear operator $T:V\rightarrow V$ such that $T(\mathcal{A})=R_{\mathcal{A}}$.

\begin{lema}\label{lemaTres}
The connectivity matrix is symmetric and is the associated matrix of the operator $T$ relative to the ordered basis $Part_{n}$.
\end{lema}
\dem
By definition, $A$ is symmetric. Relative to the ordered basis $Part_{n}$ we have $$[T(\mathcal{A}_{j})]_{Part_{n}}=[R_{\mathcal{A}_{j}}]_{Part_{n}}=(a_{1j}, a_{2j},\ldots a_{mj})$$ so $$[T]_{Part_{n}}=(a_{ij})$$
\fdem

The next proposition will be proved in a forthcoming second part paper.

\begin{prop}
The determinant of the connectivity matrix $A$ is: $$det(A)= \pm \prod_{\mathcal{A}\in Part_{n}}(m_{\mathcal{A}}-1)!$$ where $m_{\mathcal{A}}$ is the number of classes in the partition $\mathcal{A}$. In particular, $A$ is invertible.
\end{prop}

%\begin{obs}
%The operator $_{\Z}T:\Z \langle Part_{n} \rangle\rightarrow \Z \langle Part_{n} \rangle$ such that $T(\mathcal{A})=R_{\mathcal{A}}$ whose relation with the operator $T$ is $$T=\Q\otimes\ _{\Z}T$$ is not an automorphism of the $\Z$-vector space generated by $Con$ if $n>2$. In fact, its cokernel is (Proposition \ref{FactInv}) $$Coker(_{\Z}T)\simeq \frac{\Z \langle Con \rangle}{\Z \langle R_{\mathcal{A}}\ /\ \mathcal{A}\in\ Con\rangle} \ \simeq \ \bigoplus_{\mathcal{A}\in Con} \Z_{(m_{\mathcal{A}}-1)!}$$ where $m_{\mathcal{A}}$ is the number of classes in the connectivity state $\mathcal{A}$.
%\end{obs}

We observe that the left hand side of the equality is in some sense topological (it is related to connectedness) while the right hand side is combinatorial.

Finally we have the factorization theorem which solves the posed combinatorial problem. This is the main theorem of the paper.

\begin{teo}\label{MainTheorem}
Let $(b_{ij})=A^{-1}$ where $A$ is the connectivity matrix. Then $$R=\sum_{i,j=1}^{m}b_{ij}R_{\mathcal{A}_{i}} \otimes R_{\mathcal{A}_{j}}$$ and the above expression doesn't depend on the order of the basis $Part_{n}$.
\end{teo}
\dem
By Lemma \ref{lemaTres}, we have $$\mathcal{A}_{j}=\sum_{i=1}^{m} b_{ij}R_{\mathcal{A}_{i}}$$ so $$R=\sum_{h,i,j,k=1}^{m} a_{ij}b_{ki}b_{hj}R_{\mathcal{A}_{k}}\otimes R_{\mathcal{A}_{h}}$$ Because $(a_{ij})$ and $(b_{ij})$ are inverses $$\sum_{i=1}^{m} b_{ki}a_{ij}= \delta_{kj}$$ and we get $$R=\sum_{h,j,k=1}^{m} \delta_{kj} b_{hj}R_{\mathcal{A}_{k}}\otimes R_{\mathcal{A}_{h}}= \sum_{h,k=1}^{m} b_{hk}R_{\mathcal{A}_{k}}\otimes R_{\mathcal{A}_{h}}$$ Because $A$ is symmetric, its inverse is too so $$R=\sum_{h,k=1}^{m} b_{kh}R_{\mathcal{A}_{k}}\otimes R_{\mathcal{A}_{h}}$$

Finally, because the definition of $R$ doesn't depend on the order of the basis $Part_{n}$, we have the result.
\fdem

%It is interesting to write an intrinsic formulation not depending on any order of the basis:

%\begin{eqnarray*}
% \nonumber to remove numbering (before each equation)
%  R &=& \sum_{\mathcal{A}\in  Part_{n}}\mathcal{A}\otimes R_{\mathcal{A}} \\
%   &=& \sum_{\mathcal{A}\in  Part_{n}}T^{-1}(T(\mathcal{A}))\otimes R_{\mathcal{A}} \\
%   &=& \sum_{\mathcal{A}\in  Part_{n}}T^{-1}(R_{\mathcal{A}})\otimes R_{\mathcal{A}}
%\end{eqnarray*}

%Fixing an ordering in the basis, the last expression reduces to the one given in the above theorem.

As it was mentioned in the previous section, we get the solution to the probabilistic problem just applying the functional $P_{1}\otimes P_{2}$ on the last expression:

\begin{cor} Let $(b_{ij})=A^{-1}$ where $A$ is the connectivity matrix. Then $$R(G)=\sum_{i,j=1}^{m}b_{ij}R(G_{1}^{\mathcal{A}_{i}})\   R(G_{2}^{\mathcal{A}_{j}})$$ and the above expression doesn't depend on the order of the basis $Part_{n}$.
\end{cor}

The case $n=1$ is clear and reproduces the well known factorization respect to an articulation point. Let's see how the theorem works for $n=2$. Ordering the base $Part_{2}$ by $$Part_{2}= \{12, \overbrace{12}\}$$ we get the connectivity matrix $$A=\left(\begin{array}{cc}
                                                                       0 & 1 \\
                                                                       1 & 1 \\
                                                                     \end{array}
                                                                   \right)$$ and its inverse $$A^{-1}=\left(\begin{array}{cc}
                                                                                                              -1 & 1 \\
                                                                                                              1 & 0 \\
                                                                                                            \end{array}
                                                                                                          \right)$$ so
$$R=R_{\hat{12}}\otimes R_{12} + R_{12}\otimes R_{\hat{12}} - R_{12}\otimes R_{12}$$ Figure \ref{nDos} shows this factorization. See the appendix for another proof of this result. It is clear that this case is a generalization of the simple factorization factorization formula for all terminal reliability.

%\begin{figure}
%\begin{center}
%  \includegraphics[width=1.0\textwidth]{FactGraf2.png}\\
%  \end{center}
%  \caption{}\label{nDos}
%\end{figure}
Let's see the case $n=3$. Ordering the base $Part_{3}$ by: $$Part_{3}= \{123, 1\overbrace{23}, \overbrace{13}2, \overbrace{12}3, \overbrace{123}\}$$ we get the connectivity matrix $$A=\left(
                                                                \begin{array}{ccccc}
                                                                  0 & 0 & 0 & 0 & 1 \\
                                                                  0 & 0 & 1 & 1 & 1 \\
                                                                  0 & 1 & 0 & 1 & 1 \\
                                                                  0 & 1 & 1 & 0 & 1 \\
                                                                  1 & 1 & 1 & 1 & 1 \\
                                                                \end{array}
                                                              \right)$$ and its inverse

$$A^{-1}=\frac{1}{2}\left(\begin{array}{ccccc}1 & -1 & -1 & -1 & 2 \\ -1 & -1 & 1 & 1 & 0 \\ -1 & 1 & -1 & 1 & 0 \\ -1 & 1 & 1 & -1 & 0 \\ 2 & 0 & 0 & 0 & 0 \\ \end{array}\right)$$

\begin{figure}
\begin{center}
  \includegraphics[width=0.8\textwidth]{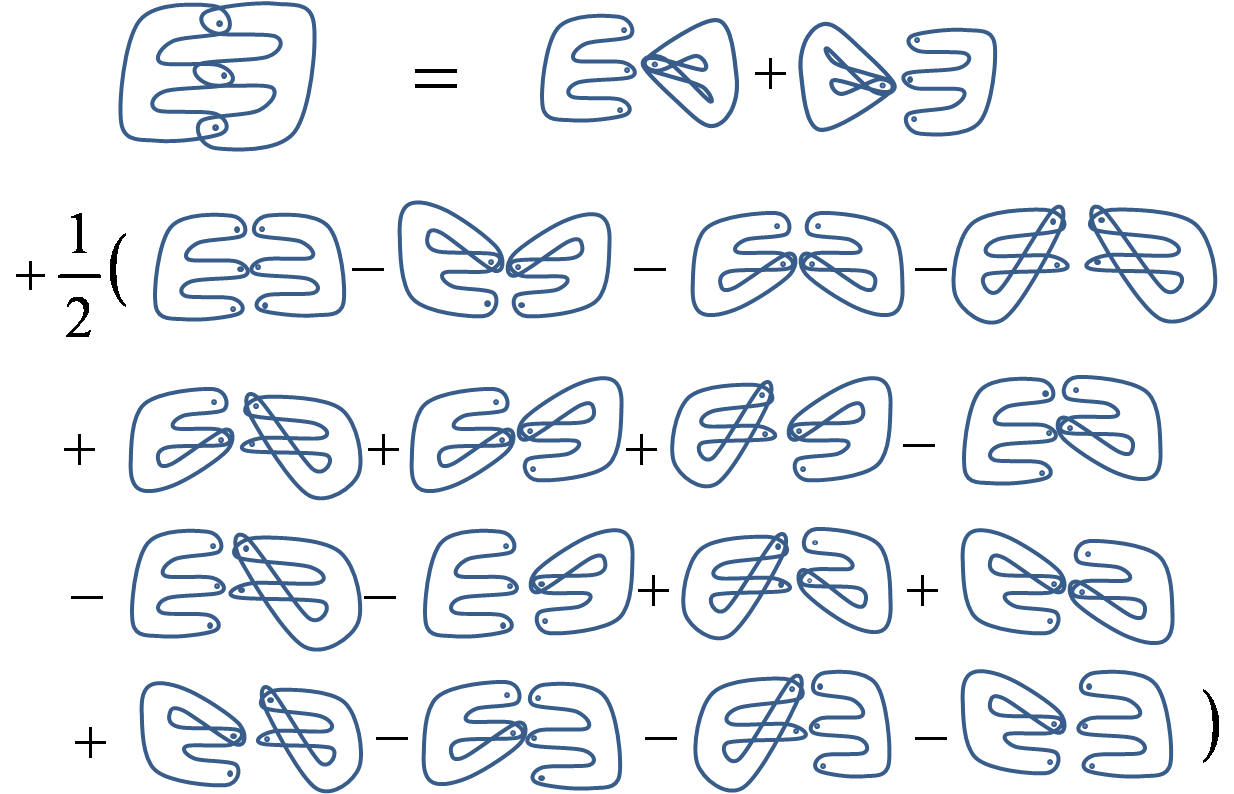}\\
  \end{center}
  \caption{}\label{nTres}
\end{figure}

The following expression for the $n=3$ factorization is illustrated in Figure \ref{nTres}:

\begin{eqnarray*}
R &=& R_{\widehat{123}}\otimes R_{123} + R_{123}\otimes R_{\widehat{123}} 
 +\frac{1}{2}(R_{123}\otimes R_{123}- R_{1\hat{23}}\otimes R_{1\hat{23}}- R_{\hat{13}2}\otimes R_{\hat{13}2}-R_{\hat{12}3}\otimes R_{\hat{12}3} \\
& & R_{1\hat{23}}\otimes R_{\hat{13}2}+ R_{1\hat{23}}\otimes R_{\hat{12}3}+ R_{\hat{13}2}\otimes R_{\hat{12}3}- R_{123}\otimes R_{1\hat{23}}- R_{123}\otimes R_{\hat{13}2}- R_{123}\otimes R_{\hat{12}3} \\
& & R_{\hat{13}2}\otimes R_{1\hat{23}}+ R_{\hat{12}3}\otimes R_{1\hat{23}}+ R_{\hat{12}3}\otimes R_{\hat{13}2}- R_{1\hat{23}}\otimes R_{123}- R_{\hat{13}2}\otimes R_{123}- R_{\hat{12}3}\otimes R_{123}) \\
\end{eqnarray*}

\section{Acknowledgements}
We are grateful to the anonymous referees for their careful reading and valuable suggestions in the improvement of the paper. The second author was supported by:
\begin{enumerate}
\item Proyecto CSIC I+D: "Propiedades de la confiabilidad di\'ametro-acotada de redes y m\'etodos de c\'alculo" \\
\item Proyect STICAmSud, tittled AMMA, institutions: INRIA Rennes-Francia, FING/UDELAR-Uruguay, UTFSM-Chile. International Coordinator, Dr. Gerardo Rubino.
\end{enumerate}

%\newpage
\appendix
\section{Another Proof of the n=2 case}

In this section, connected means edge-connected.

\begin{lema}\label{lemadefact}
Consider a stochastic graph $G= (V,E,\phi, \Phi)$ such that there exist a pair of subgraphs $G_{1}$ and $G_{2}$ of $G$ such that $E=E_{1}\cup E_{2}$, $V=V_{1}\cup V_{2}$, $E_{1}\cap E_{2}= \emptyset$ and $V_{1}\cap V_{2}= \{a,b\}$ (see the next Figure). If $G$ is connected, then $\hat{G}_{1}$ is connected and $\hat{G}_{2}$ is connected, where $\hat{G}_{i}$ is the graph resulting from the identification of the nodes $a$ and $b$ in the graph $G_{i}$.%(vea la figura (\ref{Idvert})).
\end{lema}

\begin{center}
  \includegraphics[width=0.2\textwidth]{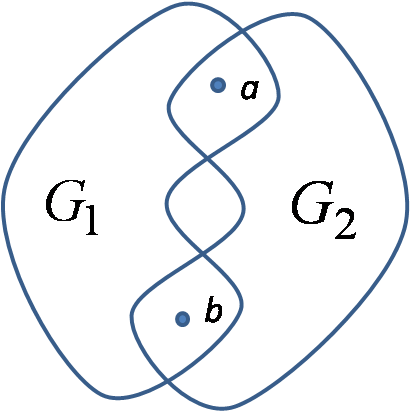}\\
\end{center}

%\begin{figure}
%\begin{center}
%  \includegraphics[width=0.6\textwidth,height=0.4\textwidth]{ContraerAristas.png}\\
%  \end{center}
%  \caption{Identificación de vértices}\label{Idvert}
%\end{figure}

\dem
Consider $v$ in $G_{1}$. There is a path in $G$ connecting $v$ with $a$. Then, there is a subpath in $G_{1}$ connecting $v$ with $a$ or $b$; i.e, there is a path in $\hat{G}_{1}$ connecting $v$ with $a=b$. This is true for every $v$ in $G_{1}$ so $\hat{G}_{1}$ is connected. In the same way as before, $\hat{G}_{2}$ is connected.
\fdem

\begin{lema}\label{lemadefactdos}
Under the hypothesis of the previous lemma, $G$ is connected if and only if one of the following items holds:
\begin{enumerate}
  \item[i.] $G_{1}$ is connected and $\hat{G}_{2}$ is connected.
  \item[ii.] $G_{2}$ is connected and $\hat{G}_{1}$ is connected.
\end{enumerate}
\end{lema}

\dem
By the previous lemma, $\hat{G}_{1}$ and $\hat{G}_{2}$ are connected. Suppose that $G_{2}$ is not connected; i.e., there are nodes $u, v$ in $G_{2}$ such that there is no path in $G_{2}$ connecting them.
\\
\\
We claim that $G_{1}$ is connected. Suppose that it is not; i.e., there are nodes $x, y$ in $G_{1}$ such that there is no path in $G_{1}$ connecting them. However, because $G$ is connectedthere exist paths $\gamma_{1}$ and $\gamma_{2}$ in $G$ connecting $x$ with $u$ and $y$ with $v$ respectively. Because $\gamma_{1}$ and $\gamma_{2}$ contain the nodes $a$ or $b$ (but not the same node), we conclude that $a$ or $b$ are not connected in $G$. In effect, if they were connected by a path in $G$, then it would exist a subpath connecting them in $G_{1}$ or $G_{2}$ so $x$ and $y$ or $u$ and $v$ would be connected. This is absurd because $G$ is connected and the nodes $a$ and $b$ are in $G$. We conclude that $G_{1}$ is connected.
\\
\\
Conversely, suppose without loss of generality that $G_{1}$ is connected and $\hat{G}_{2}$ is connected. We claim that $G$ is connected. In effect, consider $x$ and $y$ in $G$. If $x$ and $y$ are in $G_{1}$, then there is a path in $G_{1}$ connecting $x$ with $y$ and because $G_{1}\subset G$ this path is also in $G$. If $x$ is in $G_{1}$ and $y$ is in $G_{2}$, then there is a path in $G_{2}$ connecting $y$ with $a$ or $b$. Suppose without loss of generality that $y$ is connected with $a$. There is a path in $G_{1}$ connecting $a$ with $x$. The concatenation of these paths connects $x$ with $y$ in $G$. The argument for the case $x$ in $G_{2}$ and $y$ in $G_{1}$ is the same. Finally, if $x$ and $y$ are in $G_{2}$, then there are paths in $G_{2}$ connecting $x$ and $y$ with $a$ or $b$. If these paths connect $x$ and $y$ with the same point, then $x$ and $y$ are connected in $G_{2}\subset G$, otherwise there is a path in $G_{1}$ connecting $a$ with $b$ and the concatenation of these three paths connects $x$ with $y$ in $G$.
\fdem

As a corollary we have the factorization formula in the $n=2$ case (see figure \ref{nDos}).

\begin{teo}\label{teodefact}
Under the hypothesis of Lemma \ref{lemadefact} we have that $$R(G)=R(G_{1})R(\hat{G}_{2}) + R(\hat{G}_{1})R(G_{2}) - R(G_{1})R(G_{2})$$
\end{teo}
\dem
By the previous Lemma \ref{lemadefactdos}, the set of PathSets is the union of the sets of states $C_{1}$ and $C_{2}$ which verify the items $i$ and $ii$ of the lemma. The intersection of these sets is the set of states verifying that $G_{1}$ is connected and $G_{2}$ is connected. From the identity $$P(A\cup B)= P(A)+P(B)-P(A\cap B)$$ and the independence follows the result.
\fdem

As a corollary, taking $a=b$ in the previous formula we get the well known factorization respect to an articulation node.

\end{document}